\documentstyle[12pt]{article}
\newcounter{local}
\newcounter{locallocal}
\newcommand{\scl}{\stepcounter{local}}
\newcommand{\ca}{C^{\alpha, \frac{\alpha}{2}}(\bar{Q})}

\newcommand{\catau}{C^{\alpha, \frac{\alpha}{2}}(\bar{Q}_{\tau})}

\newcommand{\h}{\hspace{1cm}}
\newcommand{\hh}{\hspace{2cm}}

\newcommand{\by}{\begin{eqnarray}}
\newcommand{\ey}{\end{eqnarray}}
\newcommand{\bys}{\begin{eqnarray*}}
\newcommand{\eys}{\end{eqnarray*}}
\begin{document}

\begin{center}
{\bf On a Mathematical Model Arising from an Optimal Chemotherapeutic
Drug Treatment for Tumor Cells}

\vspace{1in}

Hong-Ming Yin \\
Department of Mathematics and Statistics\\
Washington State University\\
Pullman, WA 99164 USA.\\
Email: hyin@wsu.edu
\end{center}

\vspace{1cm} {\bf Abstract} 

In this paper we consider an optimal control problem arising from a
chemotherapeutic drug treatment for tumor cells in a living tissue. The mathematical model
for the interaction of chemotherapeutic drug and  the normal,
tumor and immune cells are governed by a nonlinear reaction-diffusion
system. We first establish the well-posedness for the nonlinear system.
Then we study the long-time behavior of the solution for the model problem. 
Finally, we design an optimal drug dosage, 
which leads to an optimal control problem under certain constrains. A complicated factor for the optimal control problem is that a minimum level of normal cells in patients must be maintained during a treatment. It is shown that there is an optimal drug dosage
for the chemotherapeutic treatment for patients.



\newpage
\begin{center}
{\bf 1. Introduction}
\end{center}

  It is well documented from the world health organization that
cancer is one of the leading causes of human morbidity and
mortality in the world. There are great efforts by scientists and medical practitioners
devoted  to find a cure  for various cancers or cancer-induced
diseases and improve the existing treatments. One of the challenging
questions for researchers is to  understand the complicated dynamical
interaction among normal, tumor, immune cells and chemotherapeutic
drugs (see \cite{EI1979,WK2014}). Toward this goal, mathematical modelling and analysis play
an important role (see \cite{F2007,F2006,RCM2007,WK2014}).

There are many mathematical models which describe cells growth
 in a living tissue. The reader can find those
Monographs by R. Rosen in 1972 (\cite{ROSEN1972}), M. Eisen in 1979 (\cite{EI1979}), V. Kuznetsov in 1992(\cite{KU1992}) and  more
recent books in 2014 by Eladdadi-Kim-Mallet (\cite{EKM2014}) and Wodarz-Komarova (\cite{WK2014}). A
well-accepted macroscopic model for cell concentration $C$ is
based on Fick's law and conservation of  mass with a logistic growth ( see \cite{ROSEN1972}):
\[ C_t=\nabla[d\nabla C]+r_0C(1-\frac{C}{K})-d_0C,\]
where $r_0$ is the coefficient of a logistic growth rate, $d$ represents the diffusion coefficient, $K$ is the maximum capacity, $d_0$ is the death rate of cell concentration.

However, when many types of cells interact each other, the model
becomes much more complicated (see \cite{RCM2007,WK2014}). In this paper
we focus only on a macroscopic interaction process for normal cell
$N$, tumor cell $T$, immune cell $I$ and a chemotherapeutic drug
$U$. Our basic goal is to understand how the chemotherapeutic drug
kills the tumor cells in an optimal way.

There are a lot of researchers studied various mathematical models about the interaction of cell dynamics
(see \cite{AD1993,BP2000,FGP1999,KP1998,OS1999,WK2014} for examples).
With this basic cell growth model in mind, several researchers developed
an ODE model to describe the interaction among normal, tumor,
immune and drug (see \cite{BP2000,KP1998,OS1999}):
\setcounter{section}{1} \setcounter{local}{1} \by
& & N_t=r_1N(1-b_2N) -c_4TN-a_3(1-e^{-U})N,\\
& & T_t=r_2T(1-b_1T)-c_2IT-c_3TN-a_2(1-e^{-U})T,\scl\\
& & I_t= s(t)+\frac{\rho IT}{\alpha+T}-c_1IT-k_1I
-a_1(1-e^{-U})I,\scl\\
& & U_t=v(t)-k_2U, \scl
\ey where $N,T,I$ and $U$ represent,
respectively, the concentration of normal, tumor, immune and
medical drug while $s(t)$ and $v(t)$ are the injection rates of
immune cell and the rate of drug.

 The values of various parameters in the ODE system (1.1)-(1.4)
 are observed and measured in medical clinics and research
 laboratories (see data in \cite{KMTP1994}).

With this ODE model, De Pillis and Radunskaya in 2003
(\cite{DR2003}) presented an very interesting dynamics about the
interaction among normal, tumor and immune cells under the
chemotherapeutic drug treatment. They obtained a range of parameters for
which that the steady-state solution of the system (1.1)-(1.4) is
stable or unstable. Particularly, they analyzed how the region of
tumor cells is changing under the influence of immune cell and
chemotherapeutic drug.

Several researchers extended the ODE models into PDE models in
order to catch the effect of diffusion process for cells in a
surrounding environment (see Ansarizadeh-Singh-Richards \cite{ASR2017},Friedman \cite{F2006},Friedman-Kim \cite{FK2011}, Lou-Ni \cite{NI},Roose-Chapman-Maini,\cite{RCM2007}, Wodarz-Komarova \cite{WK2014}). 
One of the natural PDE models in which $N,T,I$ and $U$ are governed by a
reaction-diffusion system (\cite{RCM2007}):
 \by & & N_t=\nabla
(d_1(x,t)\nabla N)+F_1(N,T,I,U), \h (x,t)\in Q, \scl\\
& & T_t=\nabla (d_2(x,t)\nabla T)+F_2(N,T,I,U), \h (x,t)\in Q,\scl \\
& &  I_t=\nabla (d_3(x,t)\nabla I)+F_3(N,T,I,U), \h (x,t)\in Q,   \scl \\
& & U_t=\nabla (d_4(x,t)\nabla U)+F_4(N,T,I,U), \h (x,t)\in Q,\scl
\ey
 where $d_i$ represents the diffusion coefficient for each type of component, $Q=\Omega\times (0, \infty)$ and $\Omega $ is a bounded domain in $R^n$.

The interaction functions are given by the same form as the ODE
model (1.1)-(1.4) (see \cite{ASR2017,CSS1992,KP1998}) with a slight modification:
\bys & & F_1(N,T,I,U)=r_1N(1-b_1N) -c_4TN-a_3(1-e^{-U})N,\\
& & F_2(N,T,I,U)=r_2T(1-b_2T)-c_2IT-c_3TN-a_2(1-e^{-U})T,\\
& & F_3(N,T,I,U)= s(x,t)+\frac{\rho IT}{\alpha+T}-c_1IT-k_1I
-a_1(1-e^{-U})I,\\
& & F_4(N,T,I,U)=v(x,t)H(N-a_0)-k_2U.
 \eys
where  various parameters in $F_i$ are derived for cell-growth
models from the clinical data (see \cite{KMTP1994}) while $H(y)$ represents the Heaviside-like  function.

Here in the PDE model we also extend the original ODE model (1.1)-(1.4) to allow the immune cells and drug inject rates depend on the location (space variables): 
\bys
 s(x,t) & = &  \mbox{growth rate of  immune
cells;}\\
v(x,t)& = & \mbox{injection rate of external chemotherapeutic
drug.} \eys

Note that $F_4$ in the current PDE model is more complicated than the original ODE model in (1.4). This is due to a clinical requirement that the drug injection must be stopped if the concentration of normal cells is lower than an acceptable level $a_0>0$.

Recently, Ansarizadeh-Singh-Richards (\cite{ASR2017}) in 2017
studied the PDE model (1.5)-(1.8) in one-space dimension without the factor $H(y$ in $F_4$. They
obtained similar dynamical stability results for the range of
various parameters as the ODE model (1.1)-(1.4). Particularly, they proved the
Jeff's phenomenon observed in clinical data. Another interesting
result in \cite{ASR2017} is that they also proved numerically that in order to slow
the growth of tumor, the chemotherapeutic drug should be injected near
the invasive front of the tumor (optimal location).

In this paper we study the model system (1.5)-(1.8) for higher space dimension subject to
appropriate initial and boundary conditions. The purpose of our
study has twofold. The first is to establish a rigorous
mathematical foundation about the well-posedness for the system (1.5)-(1.8) subject
to appropriate initial-boundary conditions. This result extended the case when the space dimension is equal to 1 (see \cite{ASR2017}). Moreover, the
global stability analysis is also established. The second is that we use some techniques developed in \cite{TR2010,YW2017} to provide a theoretical
method for an optimal drug dosage during a chemotherapeutic
tumor treatment. In the best of author's knowledge, the result is the first obtained   for the model problem with a diffusion process (see \cite{PR2001,SW1987} for the ODE case). The method could be used by medical practitioners to design an automatic drug injection device for cancer patients.

To complete the mathematical model for the nonlinear reaction-diffusion system (1.5)-(1.8), we assume that initial and boundary conditions are given by
\by && (N(x,0),T(x,0),I(x,0),U(x,0))\nonumber \\
& & =(N_0(x),T_0(x),I_0(x),U_0(x)),\,x\in
\Omega,\scl \\
& & \nabla_{\nu}(N(x,t),T(x,t),I(x,t),U(x,t))\nonumber\\
& & =(0,0,0,0),\, (x,t)\in S\times (0,\infty).\scl \ey
where $S=\partial \Omega, \nu$ is the outward unit normal on $S$ and $\nabla_{\nu}$ represents the normal derivative on $S$.

The paper is organized as follows. In section 2, we establish the well-posedness
for the nonlinear system (1.5)-(1.10). Section 3 is devoted to the stability analysis
for the solution of (1.5)-(1.10). In Section 4, we study the optimal control problem for drug dosage during a tumor treatment process. It is shown that there exists an optimal drug dosage during a chemotherapeutic treatment. Some concluding remarks are given in the final section 5.

\begin{center}
{\bf 2. Global Existence and Uniqueness}
\end{center}

In this section we show that the problem (1.5)-(1.10) admits a weak
solution globally in time under certain minimum conditions on the
known data. The definition of a weak solution for the problem (1.5)-(1.10) is standard as in \cite{Evans,LSU}. All notations for function spaces such as $L^{\infty}(Q)$ and $V_2(Q)=C([0, \infty),H^1(\Omega))$ are the same as those defined in \cite{Evans}.
For brevity, a vector function ${\bf F}=(f_1,f_2,\cdots,f_m)$ belonging in a product space $B^m$ simply means each component in the space $B$. 
The following basic conditions are assumed throughout this paper.

\ \\
H(2.1): Let $d_1(x,t), d_2(x,t), d_3(x,t)$ and $d_4(x,t)$ be of
class $L^{\infty}(Q)$ and there exist two constants $A_0$ and
$A_1$ such that
\[ 0<A_0\leq d_i(x,t)\leq A_1, \h (x,t)\in Q, i=1,\cdots,4.\]

\ \\
H(2.2): Suppose $s(x,t)$ and $ v(x,t)$ are nonnegative. There
exists a constant $A_2$ such that
\[||s||_{L^{\infty}(Q)}+ ||v||_{L^{\infty}(Q)}\leq A_2.\]
Moreover, $H(y)\in C^1(R^1)$ with $H(y)=0$ in $(-\infty,0]$ and $H(y)=1$ f in $[\delta,\infty)$ for a small $\delta>0$.

\ \\
H(2.3): Let $(N_0(x),T_0(x),I_0(x),U_0(x))\in \ca$ and there exists a constant $A_3$ such that
\[ ||N_0||_{\ca}+||T_0||_{\ca}+||I_0||_{\ca}+||U_0||_{\ca}\leq A_3.\]

 \ \\
H(2.4): All parameters $a_i, c_i, k_i, r_i$ and $\alpha, \rho$ in the system (1.5)-(1.8) are positive.

The regularity assumption in H(2.3) is not necessary for the existence of a weak solution. We give the strong regularity condition only for simplicity since our main goal is to study the optimal control problem.

In order to prove the global existence, we begin with deriving some basic a priori estimates.
In the derivation of an a priori estimate, the precise dependence of a constant $C$ is only specified at the final step.

\ \\
{\bf Lemma 2.1.} Let \[ 0<a_0\leq a(x,t)\leq a_1<\infty.\]
Let $f(x,t,u)$ be defined on $Q\times [0,\infty)$
and $f_u(x,t,u), f(x,t,0)\in L^{\infty}(Q)$ and $f(x,t,0)\geq 0$. Let $u(x,t)$ be a
 solution for the following initial-boundary value
problem: 
\setcounter{section}{2}
\setcounter{local}{1}
\by & & u_t-\nabla [a(x,t)\nabla u]
=f(x,t,u), \h (x,t)\in Q, \scl\\
& & \nabla_{\nu}u(x,t)=0, \h (x,t)\in S\times (0,\infty),\scl \\
& & u(x,0)=u_0(x)\geq 0,\h x\in \Omega. \scl \ey Suppose there exists a
constant $M$ such that
\[ \int_{\Omega} u(x,t)dx \leq M,\]
then there exists a constant $C$ such that
\[ ||u||_{L^{\infty}(Q)}\leq C(M),\]
where $C(M)$ depends only on $a, ||f_u||_{L^{\infty}(Q\times
(0,\infty))}, ||f(x,t,0)||_{L^{\infty}(Q)},
||u_0||_{L^{\infty}(\Omega)}$ and $M$.

The proof can be found in \cite{YCW2017}. A different proof can also be found in \cite{FMTY2021}. Lemma 2.1 simply indicates that a uniform estimate in $L^1(Q)$ implies a uniform estimate in $L^{\infty}(Q)$ for the solution of a linear parabolic equation with bounded coefficients.

\ \\
{\bf Lemma 2.2.} Under the assumptions H(2.1)-(2.4), a solution 
$(N,T,I,U)$ to the system (1.5)-(1.10) must satisfy the following a priori estimate:
there exists a constant $C_1$ such that
\[ 0\leq N(x,t), T(x,t),I(x,t), U(x,t)\leq C_1,\h (x,t)\in Q,\]
where $C_1$ depends only on known data.\\
{\bf Proof.}  Since we are deriving an a priori estimate, we may assume that
the solution of the system (1.5)-(1.10) is classical. 
First of all, from the nonnegativity of initial
values and $s(x,t),v(x,t)$ in $Q$, we see by the maximum principle
that
\[ N(x,t)\geq 0, T(x,t)\geq 0, I(x,t)\geq 0, U(x,t)\geq 0, \h
(x,t)\in Q.\]

 To obtain an upper bound for $U$, we note $k_2>0$ and
\[ ||v||_{L^{\infty}(Q)}\leq A_2,\]
 By applying the comparison principle with a large constant $M_1$ we see
\[ 0\leq U(x,t)\leq M_1,\]
where $M_1$ depends only on $A_1,A_2$ and
$||U_0||_{L^{\infty}(\Omega)}$.

To derive an upper bound for $N$, we  note that
\[ r_1>0, b_2>0, c_4>0, a_3>0.\]
It follows that  a large constant $M_2$ is an upper solution, where
$M_2$ depends on $r_0,a_2, b_2, c_4$ and
$||N_0||_{L^{\infty}(\Omega)}$. It follows that
\[ 0\leq N(x,t)\leq M_2, \h (x,t)\in Q.\]

Similarly, we note that $r_2>0, b_1>0, c_2>0, a_2>0$. From
\bys F_2(N,T,I,U) & = & r_1T(1-b_1T)-c_2IT-c_3TN-a_2(1-e^{-U})T\\
& \leq & r_1T(1-b_1T).
\eys
By the comparison principle, we see that a large constant $M_3$ is an upper solution.  Hence
\[ 0\leq T(x,t)\leq M_3, \h (x,t)\in Q,\]
where $M_2$ depends on $r_1, b_1, a_2, c_2$ and
$||T_0||_{L^{\infty}(\Omega)}$.

Finally, we are going to derive an upper bound for $I$. Note that
\[ c_2>0, c_4>0.\]
From the energy estimate for $I(x,t)$ in (1.6), we have
\[ \sup_{0<t<\infty}\int_{\Omega}I(x,t)dx+\int_{Q}[TI+TN]dxdt\leq
C,\] where $C$ depends only on known data.

Since $T(x,t)I(x,t)$ is uniformly bounded in $L^1(Q)$-norm, we can
use the same energy method for $I$ to obtain \bys
 & & \frac{d}{dt}\int_{\Omega} I(x,t)dx+k_1\int_{\Omega}I(x,t)dx\\
 & &  \leq  C[\int_{\Omega} s(x,t)dx+\int_{\Omega}I(x,t)T(x,t)dx]\\
 & & \leq C.
 \eys
It follows that
\[\sup_{0<t<\infty}\int_{\Omega}I(x,t)dx\leq C,\]
where $C$ is a constant which depends only on known constants.

Now we can use the result in Lemma 2.1 to obtain
\[ 0\leq I(x,t)\leq M_4,\]
where $M_4$ is a constant which depends only on known data.

\hfill Q.E.D.

With the above a priori estimate, we can establish the following global existence.

\ \\
{\bf Theorem 2.3}. Under the assumptions H(2,1)-H(2.4), the system (1.5)-(1.10) has a unique weak solution  
\[ (N,T,I,U)\in L^{\infty}(Q)\bigcap C([0,\infty); H^1(\Omega))\bigcap C^{\alpha, \frac{\alpha}{2}}(\bar{Q}). \]\\
Moreover,
\[ 0\leq N(x,t), T(x,t), I(x,t), U(x,t)\leq C,\]
where $C$ is a constant which depends only on known data. \\
{\bf Proof.} With the a priori bounds in Lemma 2.2, we can
establish the global existence by using Schauder's fixed point
theorem or a bootstrap argument (\cite{Evans}). The method is very standard. We only give an outline of the proof here.
Let $Q_{\tau}=\Omega \times (0,\tau]$ for any $\tau \in
(0,\infty)$.

Let
\bys 
 X= & & \{ (N,T,I,U)\in L^{\infty}(Q_{\tau}): N,T,I,U\geq 0,\\ 
& & ||N||_{\infty}+||T||_{\infty}+||I||_{\infty}+||U||_{\infty}\leq C_1\},
\eys
where $C_1$ is a constant in Lemma 2.2.

 It is clear that $X$ is a convex subset of
$L^{\infty}(Q_{\tau})$.

For a given $(\hat{N},\hat{T},\hat{I},\hat{U})\in X$, 
we use $(\hat{N},\hat{T},\hat{I},\hat{U})$ to replace $(N,T,I,U)$ in $F_1,F_2,F_3$ and $F_4$ and consider the corresponding linear system (1.5)-(1.10). Then from the standard theory of parabolic equations (\cite{Evans})
the corresponding linear problem
(1.5)-(1.10) has a unique solution in $C([0,\tau];H^1(\Omega))$,
denoted by $(N^*,T^*,I^*,U^*)$. Moreover, by the maximum principle,
we see
\[ 0\leq N^*(x,t),T^*(x,t),I^*(x,t),U^*(x,t)\leq C, \h (x,t)\in Q_{\tau}.\]
Furthermore, by the
DiGorgi-Nash's estimate for parabolic equations (\cite{LSU}), we
see that $(N,T,I,U)\in \ca$ and
\[||N||_{\catau}+||T||_{\catau}+||I||_{\catau}+||U||_{\catau}\leq C_0,\]
where $C_0$ depends on $ M_3$, an upper bound of $\tau$ and other
known data.

Now we define a mapping $K$ from $X$ to $\ca\subset X$ by
\[ K[(N,T,I,U)]=(N^*,T^*,I^*,U^*).\]
By Lemma 2.2, we see that all fixed points of $K$ satisfy
\[ ||N||_{\infty}+||T||_{\infty}+||I||_{\infty}+||U||_{\infty}\leq
M.\] It follows that the mapping $K$ is from $X$ into $X$.

 To show that the mapping $K$ is  continuous, we assume that
 a sequence $(N_n, T_n, I_n, U_n)\in X$ which converges to
 $(N,T,I,U)$ under the $L^{\infty}$-norm, we use
 $(N_{n}^{*},T_{n}^{*},I_{n}^{*},U_{n}^{*})$ and $(N,T,I,U)$ to
 the corresponding solutions of (1.1)-(1.6). By using standard energy estimates (\cite{LSU}), we have
 \bys
 & &
 ||N_{n}^{*}-N^*||_{\infty}+||T_{n}^{*}-T^*||_{\infty}+||I_{n}^{*}-I^*||_{\infty}+||U_{n}^{*}-U^*||_{\infty}\\
 & & \leq
C[||N_n-N||_{\infty}+||T_n-T||_{\infty}+||I_n-I||_{\infty}+||U_n-U||_{\infty}],
\eys
where $C$ depends only on known data.

 It follows that the mapping $K$ is continuous from $X$ to
$\ca$.

Finally, since the embedding operator from $\catau$ to $X$ is
compact, we see that $ K$ is a compact continuous mapping from $X$
to $X$. The Schauder's fixed point theorem yields that the mapping
$K$ has a fixed point, which is a solution of the problem
(1.5)-(1.10). Moreover, the solution $(N,T,I,U)$ satisfies
\[ ||N||_{\infty}+||T||_{\infty}+||I||_{\infty}+||U||_{\infty}\leq
M.\] Furthermore, the standard parabolic theory yields
\[||N||_{\catau}+||T||_{\catau}+||I||_{\catau}+||U||_{\catau}\leq C,\]
where $C$ depends on an upper bound of $\tau$ and other known
data.

Since $\tau$ is arbitrary, the existence of a global solution is
established. For the uniqueness, we note that
\[ \int_{\Omega}[]H(u)-H(v)](u-v))dx\leq C\int_{\Omega}(u-v)^2dx,\]
where $C$ depends only on $||u||_{L^{\infty}(\Omega)}$ and $||u||_{L^{\infty}(\Omega)}$. 

Since every weak solution is bounded,
by using the standard energy method, we see the uniqueness.

\hfill Q.E.D.

With the standard regularity theory for parabolic equations
(\cite{Evans,Lieberman,LSU}), we have the following regularity result.

\ \\
{\bf Corollary 2.4.} In addition to the conditions H(2.1)-(2.4), if all diffusion coefficients and initial values are smooth in $Q$, then the solution $(N,T,I,U)$ is classical in $Q$.

\hfill Q.E.D.

\ \\
{\bf Remark 2.1.} From the mathematical point of view, the result in Theorem 2.1 could be generalized to a more general nonlinear reaction-diffusion system. 

\begin{center}
{\bf 3. Global Stability Analysis}
\end{center}

To illustrate the basic idea,  in this section we assume that the diffusion coefficients depend only on space variables for simplicity.

\ \\
H(3.1). Let $s(x,t)$ and $v(x,t)$ converge to $v_0(x)$ and
$s_0(x)$ uniformly over $\bar{\Omega}$ as $t\rightarrow \infty$.

The following result shows that the tumor cells will be eliminated under certain conditions. 

\ \\
{\bf Theorem 3.1.} (Global Stability) Let the assumptions
H(2.1)-(2.4) and H(3.1) hold. There exists a number $r_0, \beta>0$
and such that if $0<r_2\leq r_0$ then
\bys
& & \lim_{t\rightarrow \infty}N(x,t)=N^*(x),\\
& & \lim_{t\rightarrow \infty}T(x,t)=0;\\
& & \lim_{t\rightarrow \infty}I(x,t)=I^*(x), \lim_{t\rightarrow \infty}U(x,t)=U^*(x),\h \mbox{uniformly},
\eys
where $(N^*(x), I^*(x), U^*(x))$ is the steady-state solution of the
following elliptic system:
 \setcounter{section}{3}
 \setcounter{local}{1} \by
 & & -\nabla[d_1(x)\nabla N]=F_1^*(N,U),\\
& &  -\nabla (d_3(x)\nabla I)=F_{3}^{*}(I,U), \h x\in \Omega, \scl  \\
& & -\nabla (d_4(x)\nabla U)=F_{4}^{*}(I,U), \h x\in
\Omega,\scl \\
& & \nabla_{\nu}(N, I,U)=(0,0,0),\h x\in S,\scl \ey while 
\bys 
& & F_1^*(N,I,U)=r_1N(1-b_1N)-a_3(1-e^{-U})N\\
& & F_{3}^{*}(I,U)= s_0(x)-k_1I
-a_1(1-e^{-U})I,\scl\\
& & F_{4}^{*}(I,U)=v_0(x)H(N-a_0)-k_2U.\scl
 \eys
 
 To prove Theorem 3.1, we need an elementary lemma.
 \ \\
 {\bf Lemma 3.2.} Let $d(x)\in L^{\infty}(\Omega)$ with
 \[ 0<d_0\leq d(x)\leq d_1, \h x\in \Omega\] and $k(x,t)\in L^{\infty}(Q)$
 with $k(x,t)\geq \beta>0$ in $Q$.
 Suppose $u(x,t)$ is a weak solution of the parabolic equation
 \bys
 & & u_t-\nabla (d(x)\nabla u)=-k(x,t)u, \h (x,t)\in Q,\\
 & & \nabla_{\nu}u=0, \hh (x,t)\in S\times (0,\infty),\\
 & & u(x,0)=u_0(x),\hh x\in \Omega,
 \eys
 Then
 \[ ||u(x,t)||_{L^{\infty}(\Omega)}\leq Ce^{-\beta t},\]
 where $C$ depends only on $\Omega$ and $||u||_{L^{\infty}(\Omega)}$\\
{\bf Proof.} Let $v(x,t)$ be the solution 
\bys
 & & v_t-\nabla (d(x)\nabla v)=-\beta v, \h (x,t)\in Q,\\
 & & \nabla_{\nu}v=0, \hh (x,t)\in S\times (0,\infty),\\
 & & v(x,0)=u_0(x),\hh x\in \Omega,
 \eys
Obviously (see \cite{Temam1988}), we see that $v(x,t)$ satisfies
\[||v||_{L^{\infty}(\Omega)}\leq Ce^{-\beta t}.\]
The comparison principle yields that $u(x,t)$ satisfies the same decay estimate.

\hfill Q.E.D.

\ \\
{\bf Proof of Theorem 3.1.} First of all, the maximum principle implies that
$U^*(x)$ is uniformly bounded:
\[ ||U^*||_{L^{\infty}(\Omega)}\leq C||v_0||_{L^{\infty}(\Omega)}.\]
where $C$ depends only on $\Omega$.

With a bounded $U^*(x)$, from Eq.(3.1) and Eq.(3.2) we use the comparison principle to see that
$N^*(x)$ and $I^*(x)$ are uniformly bounded in $\Omega$:
\[ ||U^*||_{L^{\infty}(\Omega)}+||I^*||_{L^{\infty}(\Omega)}\leq C,\]
where $C$ depends only on known data.

By using a similar argument to the time-dependent system, we see that the elliptic system (3.1)-(3.4) has a unique solution
$(N^*(x), I^*(x),U^*(x))\in H^1(\Omega)\bigcap C^{\alpha}(\bar{\Omega})$.
Define 
\bys
& & K(x,t):=N(x,t)-N^*(x),\\
& & J(x,t)=I(x,t)-I^*(x), \\
& & L(x,t)=U(x,t)-U^*(x), (x, t)\in Q.
\eys
Then $L(x,t)$ satisfies
\[ L_t- \nabla [d_4(x)\nabla L]=(v(x,t)H(N-a_0)-v_0(x)H(N^*-a_0))-k_2L, \h (x,t)\in
Q,\] subject to \bys & &  \nabla_{\nu} L(x,t)=0, \h (x,t)\in S\times
(0,\infty),\\
& & L(x,0)=U_0(x)-v_0(x), \h x\in \Omega. \eys
We use the maximum
principle for $W(x,t):=e^{k_{1}t}L(x,t)$ to obtain
\[ ||L||_{0}\leq C[||U_0-v_0||_0e^{-k_{1}t}+||v-v_0||_0+||N-N^*||_0],\]
where $C$ depends only on known data.

Next we note that $J(x,t)$ satisfies
\bys & & J_t-\nabla
[d_3(x)\nabla J]\\
& & =(s(x,t)-s_0(x))-k_1J-a_1[(1-e^{-U})I-(1-e^{-U^{*}})I^*]+g(x,t),\eys
where
\[ g(x,t)=\frac{\rho IT}{\alpha+T}-c_1IT.\]
Moreover,
\bys
& &  \nabla_{\nu}J(x,t)=0, \h (x,t)\in S\times (0,\infty),\\
& & J(x,0)=I_0(x)-I^*(x), \h x\in \Omega. \eys

Again we apply the maximum principle for $e^{k_{2}t}J(x,t)$ to
obtain \bys
 ||J||_{0} & \leq &
C[||I_0-I^*||_0e^{-k_{2}t}+||s-s_0||_0+||U-U^*||_0+||g||_0]\\
& \leq & C[||U_0-v_0||_0e^{-k_{1}t}+||I_0-I^*||_0e^{-k_{2}t}\\
 & & +||v-v_0||_0+||s-s_0||_0+||g||_0].
 \eys
where $||\cdot||_0$ represents the maximum norm over $\bar{Q}$.

Since $U(x,t)$ and $U^*(x)$ are bounded in $\bar{Q}$, there exists
a constant $a_0$ such that
\[ (1-e^{U}), (1-e^{-U^{*}})\geq a_0>0, \h (x,t)\in Q.\]
It follows by Lemma 3.2 that, if $ r_2$ is suitably small, we have
\[||T||_0\leq Ce^{-\beta t},\]
where $\beta>0$ is a constant which depends only on $ r_2,
a_0, a_1, a_2$ and $\Omega$.

Therefore,
\[ ||g||_{0}=o(1), \h \mbox{as $t\rightarrow \infty$}.\]
Note that $K(x,t)$ satisfies
\bys
& &  K_t-\nabla [d_1(x)\nabla K] =F_1(x,t),\h (x,t)\in Q\\
& & \nabla_{\nu}K=0, \hh (x,t)\in S\times (0,\infty).
\eys
where
\[ F_1(x,t):=f_1(N,I,T,U)-f_1^*(N^*,I^*,U^*).\]
By applying Lemma 3.2, we see  that there exists a constant $\beta_1>0$ such that
\[ ||K||_0\leq Ce^{-\beta_1 t}[||U-U^*||_0+||T||_0].\]

Finally, by H(3.1) we conclude that
\[ ||N-N^*||_0+||U-U^*||_0\leq C[||s-s_0||_0+||v-v_0||_0+o(1)], \h \mbox{as $t\rightarrow
\infty$},\]
where $C$ is a constant depending only on known data.

This concludes our desired result.

\hfill Q.E.D.
 
\begin{center}
{\bf 4. The Optimal Chemotherapeutic Drug Dosage.}
\end{center}

In this section, we introduce a method that will lead to an
optimal drug usage during a chemotherapeutic treatment.

In order to state the problem properly, we assume that
the total concentration of normal and immune cells 
maintained in a natural level when the living tissue is tumor-free.

Let
\[ A_0=\int_{\Omega}N_0(x) dx, B_0=\int_{\Omega}I_0(x)dx.\]

Let $t_0$ be a fixed time and $(N,T,I,U)$ be a solution of the system
(1.5)-(1.10).

\ \\
H(4.1): Suppose there exist constants $a_0<A_0, b_0<B_0$ such that
\[ A_0-a_0\leq \int_{0}^{t_{0}}\int_{\Omega} N(x,t) dxdt \leq A_0+a_0; B_0-b_0\leq \int_{0}^{t_{0}}\int_{\Omega}I(x,t) dxdt \leq
B_0+b_0,\] where $N$ and $I$ are the solution of the
system (1.5)-(1.10) corresponding to $U(x,t)=0$ in $ \Omega \times [0,t_0]$.

Let $Q_{t_0}=\Omega\times (0,t_0]$. We introduce an admissible set
\[ U_{ad}=\{ v(x,t)\in L^{2}(Q_{t_0}): ||v||_{L^{2}(Q_{t_0})}<\infty.\}\]

During a chemotherapeutic process, the goal is to 
find the optimal drug dosage which will minimize the total amount of tumor cells. This leads to the following optimal control problem.

For every $v\in U_{ad}$, we define the cost functional as follows:
\[ J_0(v; N,T,I,U)=\int_{\Omega}T(x,t_0)^2 dx+\frac{\lambda}{2}\int_{Q_{t_0}}v(x,t)^2dxdt,\]
where $\lambda> 0$ is a regularization parameter.

During a chemotherapeutic treatment for patients,  we have to make
sure that the normal and immune cells must maintain in an
acceptable level. Hence we impose the following  constraints:
\setcounter{section}{4} \setcounter{local}{1}
\begin{eqnarray} \int_{Q_{t_0}}N(x,t)^2 dxdt \geq A_0-a_0;\int_{Q_{t_0}}I(x,t)^2 dxdt \geq B_0-b_0,
\end{eqnarray}

Now we can state the following optimal control problem:

\ \\
{\bf Optimal Control Problem}: {\em Find $u(x,t)\in U_{ad}$ such that \by J_0 =\inf_{v(x,t)\in
U_{ad}}J(v;N,T,I,U), \scl \ey subject to the constraint (4.1),
where $(N,T,I,U)$ is a solution of the reaction-diffusion system
(1.5)-(1.10) with the drug injection $v(x,t)$.}

The main result in this section is the following theorem.

\ \\
{\bf Theorem 4.1.} Under the assumption H(4.1) and $\lambda >0$, there exists an
optimal control $u(x,t)$ for
the optimal control problem (4.2).\\
{\bf Proof.} First of all, by the assumptions for the initial data
$N_0(x)$ and $I_0(x)$, there exists a time $t_0>0$ such that
\[ \int_{0}^{t_{0}}\int_{\Omega} N(x,t)dxdt
\geq A_0- a_0, \int_{0}^{t_{0}}\int_{\Omega}
I(x,t)dxdt \geq B_0-b_0.\]

Due to the constraint (4.1), we introduce an indicator function (penal function)
$\beta(s)$:
\[
\beta(s)=\left\{ 
\begin{array}{ll}
     0,    & \mbox{if $s\geq \min\{A_0-a_0, B_0-b_0\}$ }\\
    \infty, & \mbox{otherwise}.
\end{array}
\right.
\]
The indicator function will ensure that the drug will be stopped immediately once
the amount of either normal or immune cells is below the acceptable level.

We take a smooth approximation with a small $\varepsilon >0$, denoted by $\beta_{\varepsilon}(s)$, for the indicator function $\beta(s)$ as follows:
\bys
\beta_{\varepsilon}(s)=\left\{ \begin{array}{ll}
0, & \mbox{ if $s\geq min\{A_0-a_0, B_0-b_0\}$,}\\
\mbox{smooth}, & \mbox{ if  $s< min\{A_0-a_0, B_0-b_0\}-\varepsilon $}\\
\frac{1}{\varepsilon}, & \mbox{otherwise}.
\end{array}
\right.
\eys
Note that $\beta_{\varepsilon}(s)$ implicitly depends on the lower $L^{1}(Q_{t_0})$-bounds for $U$ and $I$.

Consider the following approximate cost functional:
\[ J_{\varepsilon}(v;N,T,I,U):=J(u;N,T,I,U)+\beta_{\varepsilon}(s;U,I).\]

The cost functional $J_{\varepsilon}(v; N,T,I,U)$ is nonnegative, convex and smooth. It follows by \cite{TR2010} that there exists a minimum value $J_0(\varepsilon)$ for any $\varepsilon>0$. Set
\[ J_0(\varepsilon) =\min_{v\in U_{ad}}J_{\varepsilon}(v; N,T,I,U).\]

Suppose $u_n(x,t)$ is a minimizer sequence for $J_0$ with the underlying state solution $(N_n, T_n, I_n, U_n)$ corresponding to the system (1.5)-(1.10) with $v(x,t)=u_n(x,t)$. From the definition of the penal function $\beta_{\varepsilon}(s,U_n,I_n)$, we see $J_0(\varepsilon)\rightarrow \infty$ as long as $\varepsilon>0$ and $\varepsilon\rightarrow 0$.
It follows that the constrain (4.1) must be satisfied if $n$ is sufficiently large. Moreover, for sufficiently large $n$,
\[ J_{\varepsilon}(u_n; N_n,T_n,I_n, U_n)=J(u_n; N_n,T_n,I_n, U_n)); J_0(\varepsilon)=J_0(0).\]
Note that 
\[ ||u_n||_{L^{2}(Q_{t_0})}\leq \frac{2J_0(0)}{\lambda},\]
the weak compactness of $L^2(Q_{t_0})$ implies that there exists a limit function $u(x,t)\in L^2(Q_{t_0})$ such that 
\[ u_n(x,t)\rightarrow u(x,t),\h \mbox{weakly in $L^2(Q_{t_0})$}.\]
On the other hand, we use the energy estimate to see that  $(N_n,T_n,I_n,U_n)$ satisfies the following estimates: \bys
& & N_n(x,t), T_n(x,t), I_n(x,t), U_n(x,t)\geq 0, \h (x,t)\in Q_{t_{0}},\\
& &  ||N_n||_{V_{2}(Q_{t_0})}+||T_n||_{V_{2}(Q_{t_0})}+||I_n||_{V_{2}(Q_{t_0})}+||U_n||_{V_{2}(Q_{t_0})}\leq
C,
\eys
where $V_2(Q_{t_0})=C([0,t_0),H^1(\Omega))$ is the Sobolev space defined in $\cite{LSU}$ and $C$ is a constant depending only on known data and $\lambda$.

By the compactness embedding from $C([0,t_0],H^1(\Omega))$ into $L^2(Q_{t_0})$ (\cite{Temam1988}),  we can extract a subsequence, still denoted
by $(N_n,T_n,I_n,U_n)$, such that $(N_n,T_n,I_n,U_n)$ converges to
a limit, denoted by $(N^*,T^*,I^*,U^*)$, weakly in $V_2(Q_{\tau})$ and strongly in $L^2(Q_{t_0})$ as well as a.e. everywhere in $Q_{t_{0}}$. It follows that
\[ F_4= u_n-k_2U_n\rightarrow u(x,t)-k_2U, \mbox{weakly in $L^2(Q_{t_0})$}.\] The regularity theory for parabolic equations (\cite{Lieberman}) implies that
\[ ||U||_{\catau}\leq C,\]
for any $\tau\leq t_0.$ Again we use the comparison principle to obtain
\[0\leq  N_n,T_n, I_n\leq C,\]
where $C$ depends only on known data.

Consequently, we can use the dominate convergence theorem to obtain
\[ F_i(N_n,T_n,I_n,U_n)\rightarrow F_i(N,T,I,U) \mbox{ in strongly $L^2(Q_{t_0})$ as $n\rightarrow \infty$.}, i=1,2,3.\]

By some  routine calculations, we see that $(N^*,T^*, I^*, U^*)$ is a solution
of the system (1.5)-(1.10) corresponding to the limit control $u(x,t)$.

\hfill Q.E.D.
\ \\
{\bf Remark 4.1.} In a forthcoming paper we will study an optimal control problem in which the chemotherapeutic drug is injected through the surface of the body.

\begin{center}

{\bf 5. Conclusion }
\end{center}
In this paper we studied a nonlinear reaction-diffusion system which describes the interaction among normal, tumor, immune cells and a medical drug. The existence of a unique global solution for the nonlinear system is established by using a Schauder's fixed-point theorem. Moreover, the asymptotic behavior of the solution is obtained. Particularly, we proved that the tumor cells will be eliminated under certain conditions. We then investigated the optimal drug dosage during a chemotherapeutic process. It is proved that there exists an optimal amount of the drug dosage which will minimize the tumor cells in a fixed time period. This result could be potentially used by medical practitioners to design an automatic drug deliver device for a tumor treatment. More study with clinical data is needed for potential applications in future.

\ \\
{\bf Acknowledgement.} During the preparation of this paper, the author had several discussions with Professor Xinfu Chen at University of Pittsburgh. The author would like to express his gratitude to Professor Chen for his suggestions and comments which improved the original version of the paper.

\end{document}